\documentclass[notitlepage,leqno,10pt]{article}

\newenvironment{proof}{\noindent{\sc Proof}.\enspace}{\rule{2mm}{2mm}\medskip}

\usepackage{amsmath}
\usepackage{amssymb}

\newtheorem{theorem}{Theorem}[section]

\usepackage{amsmath}
\usepackage{amssymb}

\newcommand{\loc}{{loc}}

\newcommand{\supp}{{ \rm supp}}

\newcommand{\psido}{\psi{\rm do}}

\newtheorem{lemmy}[theorem]{Lemma}
\newcounter{proof_step}

\title{Critical regularity for   elliptic equations
from Littlewood-Paley theory II }

\author{Denis A. Labutin}

\date{}

\begin{document}
\sloppy

\maketitle

\begin{center}

{\small Department of mathematics, 
University of California, Santa Barbara,
CA 93106, USA}

\end{center}

\footnotetext[1]{E-mail address: 
labutin@math.ucsb.edu}

\

\

\noindent {\sc abstract}.
We establish a general theorem  improving regularity of solutions of   
elliptic pseudodifferential equations.
It allows to resolve in a unified way the  regularity issue
for a broad class of nonlinear elliptic equations and systems
appearing in different areas of geometry and analysis. 

\begin{center}

\bigskip\bigskip

\centerline{\bf AMS subject classification: 
35J60, 58J05, 58J40}

\end{center}

%
%
%
%
%
%
%
%
%
%
%
%
%
%
%
%

\section{Introduction}

\setcounter{equation}{0}

The purpose of this paper is to prove  a rather  general local 
regularity theorem
for elliptic equations. 
This theorem
(Theorem~\ref{maintheorem}) 
implies  the full 
$C^\infty$-regularity of solutions 
to a large class of  nonlinear elliptic problems.
For example, 
a direct consequence 
of the theorem is 
$C^\infty$-smoothness
of weak solutions  to a full scale
of nonlinear  equations with
conformally invariant pseudodifferential operators
$P^n_k$
introduced by
Paneitz, Branson, Graham, Jenne, Mason, and Sparling
\cite{Paneitz},
\cite{Branson_85},
\cite{GJMS}.
Another consequence is the critical  
regularity for elliptic systems
with the nonlinear structure analogous to the stationary 
Navier-Stokes 
system.

In the sequel  by 
$M$
we always denote a smooth  paracompact  manifold of dimension
$n$, and 
by
$X$
an open subset of
$M$.
By 
$\Psi^m(M)$
we denote the class of all pseudodifferential
operators on
$M$
of the order 
$m$,
see
section 2
for the  definitions.

%
%
%
%
%
%
%
\begin{theorem}
\label{maintheorem}
Let 
$L$,
$P$,
and
$Q$
be pseudodifferential operators from
$\Psi^\alpha(M)$,
$\Psi^\beta(M)$,
and
$\Psi^\gamma(M)$
respectively.
Assume that
$L$
is elliptic of order
$\alpha$ 
and that
\begin{equation}
\label{assumption1}
\alpha, \beta,\gamma\geq 0, \quad
\alpha>\beta+\gamma.
\end{equation}
Let 
$u$
be a distribution on
$M$
solving 
\begin{equation}
\label{maineq}
Lu
+
P(V(x)Qu)
=0
\quad
{\it on}
\quad
M.
\end{equation}
Assume that 
$V\in L^{n/(\alpha-\beta-\gamma)}(X)$
and
$u\in W^{s, p}(X)$
with
\begin{equation}
\label{technical_assumption}
\alpha-\beta\geq s\geq \gamma,
\end{equation}
and
\begin{equation}
\label{assumption2}
\gamma>s-\frac{n}{p}>\alpha -\beta -n.
\end{equation}
Then
there exists
$\varepsilon>0$,
$\varepsilon =\varepsilon(n,\alpha, \beta,\gamma, s,p)$
such that
\begin{equation}
\label{smoothness_improved}
u\in
W^{s, p+\varepsilon} (X')
\end{equation}
for any
$X'\subset\subset X$.
\end{theorem}

The assumption
$V \in L^{n/(\alpha-\beta-\gamma)}$
is the {\it critical} integrability condition.
Arbitrarily small improvement of the exponent
allows to derive
(\ref{smoothness_improved})
directly from the  Calderon-Zygmund estimate.
Simple examples show
\cite{Labutin_1}
that the integrability with any power
smaller that critical
is not enough to conclude
(\ref{smoothness_improved})
with any
$\varepsilon>0$.
According to the general principle,
results
for linear problems 
of type 
(\ref{maineq})
with rough coefficients
have immediate 
applications to the 
corresponding nonlinear
problems
\begin{equation*}
Lu+f(\partial, u) =0.
\end{equation*}
This paper is a continuation of our earlier work
\cite{Labutin_1}.

The main point of Theorem~\ref{maintheorem}
is that the regularity improvement holds
for a
critical equations  
{\it 
without any 
structure assumptions}
on the nonlinearity 
provided 
the nonlinearity satisfies
\begin{equation}
\label{main_point}
VQ(u)\in L^r_\loc
\quad
{\rm 
with
}
\quad
\infty>r>1.
\end{equation}
As we discuss below, several 
nonlinear  problems satisfy this condition.
Main requirement 
(\ref{assumption2})
in Theorem~\ref{maintheorem}
is nothing but 
(\ref{main_point}).
Condition
(\ref{technical_assumption})
is not important and always holds in nonlinear applications.

The question of local regularity  for critical elliptic
equations and systems is important for  
different areas of geometry and analysis.
Many works are dedicated to this problem
for various equations. The critical  equations which 
attracted a lot of 
attention include,
for example,
Yamabe equation, the family of $Q$-curvature equations,
$H$-surfaces system,
system of harmonic maps from the plane,
elliptic Yang-Mills system in dimension $4$,
and others.

The proof Theorem~\ref{maintheorem}
works also for   $\psido$
between  (complex)
vector bundles
over
$M$.
For simplicity
we do not formulate the theorem in 
the most general form.
Below we discuss classes of nonlinear  problems 
covered and not covered by the theorem.

Here is a typical scheme
for nonlinear applications of Theorem~\ref{maintheorem}.
Consider, for example,  
the following  equation
in
$\mathbf{R}^4$
\begin{equation*}
\Delta^2 u + \partial_1^2
\left((\partial_1 u)^2 \right)=0,
\quad
u\in H^2_\loc.
\end{equation*}
It 
is critical for
$H^2$.
Setting
$V=\partial_1 u$
we arrive, due to the Sobolev embedding, 
at
\begin{equation*}
\Delta^2 u +\partial_1^2
(V\partial_1 u)=0, 
\quad
u\in H^2_\loc,
\quad
V\in L^4_\loc.
\end{equation*}
Apply Theorem~\ref{maintheorem}
to 
derive that
$u\in W^{2, 2+\varepsilon}_\loc$.
Thus we improved the integrability of
$\partial_1 u$
above the critical.
Now repeated applications of Calderon-Zygmund and Schauder estimates
to the nonlinear problem imply that
$u\in C^\infty_\loc$.

As a direct consequence of
Theorem~\ref{maintheorem}
we recover the critical regularity of 
stationary Navier-Stokes system established in
\cite{vonWahl}, 
\cite{Galdi}.
That is, take any
$u=(u_1, \ldots, u_n)$,
$u_j\in H^{n/4}_\loc$,
solving
\begin{equation*}
\Delta u -\mathbf{P}\big( u\cdot\nabla\, u\big) =0,
\quad
\nabla\cdot u =0.
\end{equation*}
Here
$\mathbf{P}$
is the Leray projector on the space of divergence free vector fields. 
Operator
$\mathbf{P}$
is  a homogneous
Fourier multiplier of order 
$0$
with the singularity at the origin consistent
with the scaling.
Such  operators enjoy the same mappping properties
as
$\psido$
from
$\Psi^{0}$,
cf.
\cite{Stein}, Ch. 6.
Setting, say,
$V=u$ and
applying
Theorem~\ref{maintheorem}
deduce that
$u\in W^{n/4, 2+\varepsilon}_\loc$.
Then the application of the standard Calderon-Zygmund and 
Schauder estimates
gives
$u\in C^\infty_\loc$.
In particular,
every finite energy solution
$u\in H^1_\loc$
is smooth in dimension
$4$.

Another consequence
of Theorem~\ref{maintheorem}
is the full regularity for a class of 
nonlinear equations arising in conformal geometry.
Locally these equations have the 
$\psido$
$(-\Delta)^{k/2}$,
$k=1$,
$2$,
$\ldots$,
as the leading part.
Study of these equations and related conformal concepts
(for example, Branson's $Q$-curvature)
is an active area of current research
in
geometry, nonlinear PDEs, mathematical physics, and spectral
theory, see e.g.
\cite{Alice_Chang_Paul_Yang_ICM}
\cite{Alice_Chang_book},
\cite{Fefferman_Graham_1},
\cite{Fefferman_Graham_2},
\cite{Graham_Zworski}.
A corollary of Theorem~\ref{maintheorem}
is the local regularity
of weak solutions to the corresponding critical
nonlinear equations. For example, any 
$u\in H^{n/2}_\loc$
solving
\begin{equation*}
\Delta^{n/2} u  
+
\mathrm{div}
\left(
\left(
\Delta^{3/2}
u
\right)^{(n-2)/3}
\mathrm{grad} u
\right)
=0
\end{equation*}
is smooth (say, in dimensions $5$,
$8$,
$11$, 
...).
For {\it differential}
operators this was established in
\cite{Chang_Gursky_Yang}
(for minimising solutions),
\cite{Uhlenbeck_Viaclovsky},
\cite{Djadli_Hebey_Ledoux}
(general case).
The last two works  use  the same idea
which is based on  the unique solvability 
of the suitable elliptic  boundary value problem.
In the case of
$\psido$
this construction does not work. 
This provided the initial motivation for 
our paper.
Paper
\cite{YYLi}
contains an approach to the local regularity
of 
(\ref{maineq})
based on rewriting the equation
as an integral equation in the physical space.
Such approach works particularly well
in the case when
$L$,
$P$,
and
$Q$
have some structure
(say,  the lack of
$x$-dependence, or possession of  
a radial fundamental solution, or homogeneity, ...).
In the case of the second order single equations
DeGiorgi-Moser tools are available
\cite{Brezis_Kato},
\cite{Trudinger}.

What  critical problems are not
covered by Theorem~\ref{maintheorem}?
The examples are the equation for 
weakly harmonic maps in the plane, and the equation for  
$H$-surfaces.
The  weakly harmonic maps from 
the plane into spheres are solutions 
$u\in H^1(\Omega, \mathbf{S}^N)$,
$\Omega\subset\mathbf{R}^2$,
of the system
\begin{equation*}
\Delta u + |\nabla u|^2 u =0.
\end{equation*}
Solutions to the  
$H$-surfaces 
equation are the  functions
$u\in H^1(\Omega, \mathbf{R}^3)$
satisfying
\begin{equation*}
\Delta u + \partial_x u\wedge \partial_y u =0.
\end{equation*}
For both equations, 
solutions of the class
$H^1_\loc$
do not satisfy 
the main assumption
(\ref{assumption2})
of Theorem~\ref{maintheorem}.
For example, if we utilise the form
(\ref{maineq})
with
$P=1$,
then
the nonlinearity is only in
$L^1_\loc$.
If, alternatively, we utilise the divergence structure
of the nonlinearity to arrive at
(\ref{maineq})
with
$P\in \Psi^1$,
then the first inequality in
(\ref{assumption2})
is violated.
Therefore, one needs  to exploit more
information 
about the nonlinearity.
Presently several different proofs of the  full regularity 
for  both equations are known
\cite{Struwe_book_1},
\cite{Struwe_book_2},
\cite{Helein_book},
\cite{Tao}
\cite{Chang_Yang_Wang},
\cite{Streletzki}.
The common element in  all of them
is the heavier utilisation
of the div-curl structure
of the nonlinearity.

\section{Preliminaries}
\label{prelims}
\setcounter{equation}{0}

In this section we review some definitions and 
state some  estimates needed  in the proof of the theorem.
We refer to monographs
\cite{Sogge},
\cite{Stein},
and
\cite{Taylor}
for more information.

\subsection{Littlewood-Paley decomposition}

Let
$\{\widehat{\varphi}_j\}_{j=-\infty}^{+\infty}$
be the standard 
smooth partition of unity in the 
Littlewood-Paley theory
\cite{Sogge},
\cite{Stein}.
Thus 
$\widehat{\varphi}_j=\widehat{\varphi}(\cdot/ 2^j)$
is supported in, say,
the ring
\begin{equation*}
\{ \xi\colon 2^j 3/5 \leq |\xi| \leq 2^j 5/3\}
\subset
\left(
B_{2^{j+1}}\setminus B_{2^{j-1}}
\right).
\end{equation*}
Let 
$P_j$
denote
the Littlewood-Paley projection,
\begin{equation*}
(P_j f)^\wedge
=\widehat{\varphi}_j \widehat{f},
\quad
f\in
\mathcal{S}'.
\end{equation*}
We also set
\begin{equation*}
P_{a<\cdot<b}
=
\sum_{j=a+1}^{b-1}
P_j.
\end{equation*}
It follows easily that for any
$j$
and any
$p$,
$1\leq p\leq\infty$, 
the operators
\begin{equation*}
P_j\colon L^p\longrightarrow L^p
\end{equation*}
have norms bounded uniformly
over
$j$.

Distributions  with the localised Fourier transform 
(for example any $P_j f$)
enjoy the important {\it Bernstein inequality}.
It says that
for
$f\in\mathcal{S}'$
and
$1\leq p\leq q\leq\infty$
\begin{equation*}
\|f\|_q
\lesssim
2^{nj((1/p)-(1/q))}
\|f\|_p
\quad
{\rm provided}
\quad
\supp\,\widehat{f}
\subset
B_{2^j}.
\end{equation*}

For
$s\in\mathbf{R}^1$
and
$1<p<\infty$
the Sobolev space
$W^{s,p}=W^{s,p}(\mathbf{R}^n)$
consists of all distributions with the finite norm
defined by
\begin{equation*}
\|f\|_{W^{s,p}}^p
=
\|P_{\cdot\leq 0} f\|_p^p
+
\left\|
\sum_{j=1}^\infty
\left(
2^{2js}
|P_j f|^2
\right)^{1/2}
\right\|_p^p.
\end{equation*}
It is well known in analysis that this norm
is difficult to work with.
However,
in this paper we will use only  the following facts about
$W^{s,p}$
which follow easily
from the properties of
$P_j$.
Namely, 
\begin{eqnarray*}
&
&
f\in W^{s,p}
\Longrightarrow
\|P_k f\|_p\lesssim \frac{\|f\|_{W^{s,p}}}{2^{sk}}
\
{\rm for \ all \ }
k\geq 0;
\\
&
&
\|P_k f\|_p\leq \frac{N}{2^{(s+\varepsilon)k}},
\
\varepsilon>0,
\
{\rm for\  all}
\ 
k\geq 0
\Longrightarrow
\|f\|_{W^{s,p}}
\leq C_{\varepsilon} N.
\end{eqnarray*}
The central result of the Littlewood-Paley theory 
implies that for
$s=0$,
$1$,
$\ldots$
the space
$W^{s,p}$
consists of all  distributions  with all derivatives
up to the order
$s$
lying
in
$L^p$,
$1<p<\infty$.
For 
$p=2$
this follows from the Plansherel isometry.

If
$M$
is a smooth paracompact manifold of dimension
$n$, 
then
the Sobolev spaces
$W^{s,p}(M)$
(and other spaces of functions)
are defined via the partition of unity.
If the functions are defined on
$\mathbf{R}^n$
then  frequently we  will  not  write the domain of the definition.

\subsection{Pseudodifferential operators}

A pseudodifferential operator
($\psido$ for short)
$A$
of order
$m\in\mathbf{R}$
on
$M$
is a linear map
\begin{equation*}
A\colon C^\infty_0(M)\longrightarrow C^\infty_\loc(M)
\end{equation*}
such that in any coordinate chart
$A$
is a $\psido$
of order
$m$
in
$\mathbf{R}^n$.
This means that if we take a
coordinate patch
\begin{equation*}
U\subset M,
\quad
\kappa\colon U \to {\Omega},
\quad
{\Omega}\subset\mathbf{R}^n.
\end{equation*}
and fix
two functions
$\varphi,\psi\in C^\infty_0(U)$,
such that
$\varphi =1$
on
$\supp \psi$,
then  for
$v\in C^\infty_0 (\mathbf{R}^n)$
the
map
\begin{equation}
\label{defpsido}
v\longmapsto
\Big(
\psi A(\varphi v\circ\kappa)\Big)(\kappa^{-1} x),
\quad 
x\in {\Omega},
\end{equation}
is given by the formula
\begin{equation*}
v
\longmapsto
\int_{\mathbf{R}^n}
e^{ix\xi}
a(x,\xi)
\hat{v}(\xi)
\,
d\xi,
\end{equation*}
where
$a$
is a symbol from
$S^m$.
The latter means that
$a\in C^\infty_\loc(\mathbf{R}^n\times\mathbf{R}^n)$
and the estimate
\begin{equation*}
\left|
\partial_\xi^k\partial_x^l
(x,\xi)
\right|
\leq
C_{k,l}
(1+|\xi|)^{m-|k|}
\end{equation*}
holds 
for all
multiindicies
$k$
and
$l$.
From
(\ref{defpsido})
it follows that
$\supp a\subset\kappa(\supp \psi)\times\mathbf{R}^n$.
By
$\Psi^m(M)$
we denote the set of all
$\psido$
on
$M$
of order
(at most)
$m$.
Throughout the paper we abuse notations by denoting the operator
$A_{\kappa, \varphi, \psi}$
in
(\ref{defpsido})
by the same letter
${A}$.

We say that
$A\in \Psi^m(M)$
is elliptic of
order
$m>0$
if in  any coordinate patch
as above,
for any
$\Omega'\subset\subset\Omega$
there are constants
$C_{1,2}>0$
the  symbol of $A$
satisfies
\begin{equation*}
|a(x,\xi)|\geq C_1|\xi|^m
\quad
{\rm if}
\quad
x\in\Omega', 
\quad
|\xi|\geq C_2.
\end{equation*}

Any
$A\in \Psi^m(M)$,
$m\in\mathbf{R}$,
admits a unique extension as
\begin{equation*}
A\colon \mathcal{E}'(M)\longrightarrow \mathcal{D}'(M).
\end{equation*}
Moreover, the central result of the Calderon-Zygmund theory
can be stated as the following mapping property of
$\psido$:
\begin{equation*}
A\colon W^{s,p}(M)
\longrightarrow
W^{s-m, p}_\loc(M)
\quad
{\rm for
\quad
}
s\in\mathbf{R},
\quad
1<p<\infty.
\end{equation*}

\subsection{Estimates for pseudodifferential operators}

It is easy to see that 
$P_k\in \Psi^{-\infty}$.
For any 
$A\in\Psi^m$
and any
$k$
the estimate
\begin{equation}
\label{AP_j_estimate}
\|A P_k f\|_p 
\leq 
C_{A,p} 2^{km} 
\| P_{k-1\leq\cdot\leq k+1} f\|_p,
\quad
f\in\mathcal{S}',
\end{equation}
holds for
any
$p$,
$1\leq p\leq \infty$,
with the constant
$C_{A,p}$
independent of 
$k$.
Estimate
(\ref{AP_j_estimate})
is a consequence of 
the Calderon-Zygmund
theory
\cite{Stein}, Ch. 6.

\begin{lemmy}
\label{commutator_lemma}
For any 
$\psi{\it do}$
$A\in\Psi^m$,
and any
$P_k$,
$k\geq 10$
the estimate
\begin{eqnarray}
\label{commutator_estimate}
\|(P_k A - A P_k)f\|_p
\leq
C_{A,N}
\Big(
&
&
2^{k(m-1)}\| P_{k-5<\cdot<k+5} f\|_p
\nonumber
\\
&
&
+
2^{-Nk}
\|P_{\cdot\leq 0} f\|
\nonumber
\\
&
&
+
2^{-Nk}
\sum_{j=1}^\infty
2^{-Nj}
\|P_j f\|_p
\Big), 
\end{eqnarray} 
holds for 
$f\in\mathcal{S}'$,
and all
$N>0$,
$1\leq p\leq\infty$.
The constant
$C_{A,N}$
does not depend on
$k$.
\end{lemmy}

It is easy to see, that
despite the symbol
$\varphi_k \in S^{-\infty}$,
it admits the pointwise estimates {\it 
independent of}
$k$
only in
$S^0$.
It is due to this fact  
that we can improve only by the factor
$2^{-k}$
in
(\ref{commutator_estimate}).

\setcounter{proof_step}{1}
\begin{proof}
{\bf \arabic{proof_step}.}
\stepcounter{proof_step}
Let
$a$
be the symbol of
$A$,
$a\in S^m$.
Let
$c$
be the symbol of
$P_k A$.
The basic theorem of
$\psido$-calculus 
asserts that
$c\in S^{-\infty}$, 
and  moreover
\begin{equation*}
c(x,\xi) = 
\varphi(2^{-k}\xi)
a(x,\xi)
+
\rho_k(x,\xi)
\end{equation*}
with
$\rho_k\in S^{-\infty}$.
However, we are interested in
explicit dependence on
$k$.
We claim that
for any
multiindices
$\alpha$
and
$\beta$, 
$N>0$,
and all
$k\geq 10$
\begin{eqnarray}
\left|
\partial^\alpha_\xi \partial^\beta_x
\rho_k (x,\xi)
\right|
&
\leq
&
C_{N, A, \alpha, \beta}
\nonumber
\\
&
&
\times
\left\{
\begin{array}{lll}
(1+|\xi|)^{m-1-|\alpha|}
&
{\rm if}
&
2^{k-3}\leq |\xi| \leq 2^{k+3}
\\
(1+|\xi|)^{-N}
&
{\rm if}
& 
|\xi|\geq 2^{k+3}
\\
2^{-Nk}
&
{\rm if}
&
|\xi|\leq 2^{k-3}.
\end{array}
\right.
\label{P_jA_reminder_symbol_estimate}
\end{eqnarray}
The point here is that
$C_{N, A,\alpha, \beta}$
is independent of
$k$.

It is obvious that estimate
(\ref{commutator_estimate})
follows from
(\ref{P_jA_reminder_symbol_estimate})
and 
(\ref{AP_j_estimate}).
We will prove
(\ref{P_jA_reminder_symbol_estimate})
only for
$\alpha=\beta=0$.
The obvious changes make the proof to work for the general case.

{\bf \arabic{proof_step}.}
\stepcounter{proof_step}
To show
(\ref{P_jA_reminder_symbol_estimate})
we first write down
$c(x,\xi)$.
We can assume that
$a(x,\xi)$
has a  compact support in
$x$. 
A well-known argument allows to remove this assumption 
while preserving
(\ref{P_jA_reminder_symbol_estimate}),
see
\cite{Stein}, Ch. 6.
We also localise 
$P_k$
in 
$x$
by considering
\begin{equation*}
\varphi_{k,\varepsilon}(x,\xi)
=
\theta(\varepsilon x)\varphi(2^{-k}\xi),
\end{equation*}
where 
$\theta$
is a cut-off and
$\varepsilon>0$.
All our estimates will be independent of
$\varepsilon$,
$0<\varepsilon<1$.
The 
Fourier transform easily
gives that
\begin{eqnarray*}
c_{\varepsilon}(x,\xi)
&
=
&
\frac{1}{(2\pi)^n}
\int_{\mathbf{R}^n}
e^{i x\eta}
\varphi_{k,\varepsilon}(x, \xi+\eta)
\hat{a}(\eta, \xi)
\,
d\eta
\\
&
=
&
\varphi_{k,\varepsilon}(x, \xi)
a(x, \xi)
\\
&
&
+
\frac{1}{(2\pi)^n}
\int_{\mathbf{R}^n}
e^{i x\eta}
\hat{a}(\eta, \xi)
\left(
\varphi_{k,\varepsilon}(x, \xi+\eta)
-
\varphi_{k,\varepsilon}(x, \xi)
\right)
\,
d\eta,
\end{eqnarray*}
where
$\hat{a}(\cdot,\cdot)$
denotes the Fourier transform in  the 
{\it first} 
variable.
The localisation of
$a$
in the first variable imply that
\begin{equation}
\label{estimate_for_a}
\left|
\hat{a}(\eta, \xi)
\right|
\leq C_{A,N}
\frac{(1+|\xi|)^m}{(1+|\eta|)^N},
\quad
{\rm for\quad all\quad  }
N.
\end{equation}
Set
\begin{equation*}
R_{k,\varepsilon}(x,\xi,\eta)
=
\varphi_{k,\varepsilon}(x, \xi+\eta)
-
\varphi_{k,\varepsilon}(x, \xi)
\end{equation*}
for the reminder, which we need to estimate.

{\bf \arabic{proof_step}.}
\stepcounter{proof_step}
Let us show the first estimate in
(\ref{P_jA_reminder_symbol_estimate}).
Fix  any
$\xi$
satisfying
${2^{k-3}\leq |\xi| \leq 2^{k+3}}$.
Assume first that
$|\eta|\lesssim |\xi|$.
Then we have 
\begin{eqnarray*}
\left|
R_{k,\varepsilon}
(x,\xi,\eta)
\right|
&
\lesssim
&
\frac{|\eta|}{2^k}
\\
&
\lesssim
&
\frac{|\eta|}{1+|\xi|}.
\end{eqnarray*}
Assume next that
$|\eta|\gtrsim |\xi|$.
Then due to the localisation of
$\varphi_k$
we have
\begin{eqnarray*}
\left|
R_{k,\varepsilon}
(x,\xi,\eta)
\right|
&
=
&
|\varphi_{k,\varepsilon}(x,\xi)|
\\
&
\lesssim
&
1.
\end{eqnarray*}
Combining with
(\ref{estimate_for_a})
discover that
\begin{eqnarray*}
|\rho_{k,\varepsilon}(x,\xi)|
&
\lesssim
&
C_{A,N}
(1+|\xi|)^{m-1}
\int_{\{|\eta|\lesssim |\xi|\}}
\frac{1}{(1+|\eta|)^N}
\,
d\eta
\\
&
&
+
C_{A,N}
\int_{\{|\eta|\gtrsim |\xi|\}}
\frac{(1+|\xi|)^{m} }{(1+|\eta|)^N}
\,
d\eta
\\
&
\lesssim
&
C_{A,N} (1+|\xi|)^{m-1}.
\end{eqnarray*}
Finally let
$\varepsilon\to 0$.
Thus
the first part of 
(\ref{P_jA_reminder_symbol_estimate})
is proved.

{\bf \arabic{proof_step}.}
\stepcounter{proof_step}
Let us show the second  estimate in
(\ref{P_jA_reminder_symbol_estimate}).
Fix 
$\xi$,
$|\xi | \geq 2^{k+3}$
Then 
\begin{equation*}
\left|
R_{k,\varepsilon}
(x,\xi,\eta)
\right|
=
0
\end{equation*}
for all
$\eta$
such that
$|\eta|\lesssim |\xi|$.
Assume next that
$|\eta|\gtrsim |\xi|$.
Then
\begin{equation*}
\left|
R_{k,\varepsilon}
(x,\xi,\eta)
\right|
\lesssim
1.
\end{equation*}
Combining with
(\ref{estimate_for_a})
derive that
\begin{eqnarray*}
|\rho_{k,\varepsilon}(x,\xi)|
&
\lesssim
&
C_{A,N}
(1+|\xi|)^{m-1}
\int_{\{|\eta|\gtrsim |\xi|\}}
\frac{1}{(1+|\eta|)^N}
\,
d\eta
\\
&
\lesssim
&
C_{A,N} (1+|\xi|)^{-N'}
\
{\rm for \ any\ } N'.
\end{eqnarray*}
After letting
$\varepsilon \to 0$
we conclude that
the second part of 
(\ref{P_jA_reminder_symbol_estimate})
is also proved.

{\bf \arabic{proof_step}.}
\stepcounter{proof_step}
To show the third estimate in
(\ref{P_jA_reminder_symbol_estimate})
it is enough to notice that
for any fixed
$\xi$,
$|\xi|\leq 2^{k-3} $,
we have
\begin{equation*}
\left|
R_{k,\varepsilon}
(x,\xi,\eta)
\right|\ne 0
\quad
{\rm 
only 
\quad 
for
\quad
}
|\eta|\asymp 2^k.
\end{equation*}
Hence
\begin{eqnarray*}
|\rho_{k,\varepsilon}(x,\xi)|
&
\lesssim
&
C_{A,N}
\int_{\{|\eta|\asymp 2^k\}}
\frac{(1+|\xi|^m)}{(1+|\eta|)^N}
\,
d\eta
\\
&
\lesssim
&
C_{A,N'} 2^{-N'k}
\ \
{\rm for \ \ any\ \  } N'.
\end{eqnarray*}
Finally we let
$\varepsilon \to 0$.
\end{proof}

As a direct  consequence of 
(\ref{commutator_estimate})
we obtain that under the assumptions
of Lemma~\ref{commutator_lemma}
the estimate
\begin{equation}
\label{commutator_estimate_rough}
\|P_k A f\|_p
\leq
C_{A,N}
\left(
2^{mk}\| P_{k-10\leq\cdot\leq k+10} f\|_p
+
2^{-Nk}
\|f\|_p
\right)
\end{equation}
holds.
Estimates
(\ref{AP_j_estimate})
and
(\ref{commutator_estimate_rough})
show that 
$A\in \Psi ^m$
acts on the part of
$f$
with the frequencies
of the order
$2^k$,
as a multiplication by
$2^{mk}$
modulo
a fixed frequency spreading and an arbitrary small correction.

\section{Proof of Theorem~\ref{maintheorem}}
\setcounter{equation}{0}

%
%

\setcounter{proof_step}{1}

By
$q$
we always denote the critical exponent,
\begin{equation*}
q
=
n/(\alpha -\beta -\gamma)
.
\end{equation*}
According to
(\ref{assumption1})
and
(\ref{assumption2})
we have
$1<q<\infty$.
In this section we write
$A\lesssim B$
if
$A\leq CB$
with a  constant 
$C>0$,
depending on the operators
$L$,
$P$,
$Q$, and
the parameters from Theorem~\ref{maintheorem}.

\begin{proof}[of Theorem~\ref{maintheorem}]
{\bf \arabic{proof_step}.} 
\stepcounter{proof_step}
When proving
(\ref{smoothness_improved})
we can can assume that
$X$
lies in a single coordinate chart.
First we make  a suitable localisation of equation
(\ref{maineq}).
This will be done in several steps exploiting the pseudolocal 
character of $\psido$ in a standard way.

Apply definition
(\ref{defpsido})
and
utilise the pseudolocality 
of 
$\psido$ 
to derive that
\begin{eqnarray*}
\psi L(\varphi u) +\psi P(V(x) Qu)
&
=
& 
\psi L((\varphi -1)u)
\\
&
=
&
\psi f,
\quad
\psi f\in C^\infty_0 (\Omega).
\end{eqnarray*}
Abusing notations
we  have denoted
$u\circ\kappa^{-1}$
by the same letter 
$u$.
Repeat this argument several times to find
functions
$\psi_j\in C^\infty_0(\Omega)$,
$j=1,\ldots, 4$,
such that
$\psi_{j}=1$
on
$\supp \psi_{j+1}$,
and
\begin{equation}
\label{localised_1}
\psi_4 L(\psi_1 \, u)
+
\psi_4 P
\left(
\psi_3
\,
V(x)
\psi_2
Q(
\psi_1\, u
)
\right)
=
\psi_4 P(\psi_3\, v)
+\psi_4 f,
\end{equation}
where
\begin{equation*}
v\in L^{n/(\alpha-\beta-\gamma)}(\supp \psi_3),
\quad
f\in C^\infty_0(\supp \psi_4).
\end{equation*}

Next, 
we can assume that
$B_1\subset\subset\Omega$
and
\begin{equation*}
\psi_4 =1 \quad
{\rm on}
\quad
B_1.
\end{equation*}
Take a cutoff function
$\eta_\rho$,
\begin{equation*}
\eta_\rho =1
\quad
{\rm in}
\quad
B_\rho,
\quad
\eta_\rho =0
\quad
{\rm outside}
\quad
B_{2\rho}.
\end{equation*}
Later we will choose
$\rho$
small.
The commutator of
the multiplication by 
$\eta_{\rho}$
and 
an operator with the symbol from
$S^m$
is a 
$\psido$ 
with the symbol in
$S^{m -1}$.
For differential operators 
this is just the Leibnitz rule.
Fix a function
$\phi\in C^{\infty}_0(\Omega)$
such that
\begin{equation*}
\phi =1
\quad
{\rm on}
\quad
\supp \psi_1.
\end{equation*}
Then  we can continue
(\ref{localised_1})
and 
write
\begin{equation}
\label{interm_localisation}
L(\eta_\rho u) 
+
\eta_\rho
{P}
\left(
V(x) {Q}(\phi u)
\right)
=
\Phi_{\alpha -1}(\phi u)
+
\Phi_\beta v
+
f.
\end{equation}
From now on we denote by
$\Phi_m$
a generic operator
from
$\Psi^m$
compactly 
$x$-supported in
$\Omega$.
In 
(\ref{interm_localisation})
we have
$f\in C^\infty_0(\Omega)$, 
$v\in L^{q}$,
$\supp v\subset\subset\Omega$.

Repeat this argument several times   
utilising 
$\eta_\rho =
\eta_\rho
\eta_{2\rho}$.
Thus we
discover 
the following relation in
$\mathbf{R}^n$:
\begin{eqnarray}
{L}(\eta_\rho u) 
+
{P}
\left(
\eta_{2\rho}
V(x) 
{Q}(\eta_\rho u)
\right)
=
\Phi_{\alpha -1}(\phi u)
&
+
&
\Phi_{\beta -1} 
\left(
V(x) {Q}(\phi u)
\right)
\nonumber
\\
&
+
&
{P}
\left(
\eta_{2\rho}
V(x)
\Phi_{\gamma-1}
(\phi u)
\right)
\nonumber
\\
&
+
&
\Phi_\beta v
\nonumber
\\
&
+
&
f,
\label{localised_2}
\end{eqnarray}
where
the functions
$v$ 
and
$f$
are the same as in
(\ref{localised_1}),
(\ref{interm_localisation}).

{\bf \arabic{proof_step}.} 
\stepcounter{proof_step}
Next we utilise the ellipticity to 
rewrite the equation
in an essentially  invertible form.
The
ellipticity of 
$L$
allows us to split
it,
\begin{equation*}
{L} = E+M_\alpha,
\end{equation*}
into an operator
$E$
with the symbol
$e(x,\xi)\in S^\alpha$,
\begin{equation*}
|e(x,\xi)|\gtrsim (1+|\xi|)^\alpha
\quad
{\rm 
for
\quad 
all 
}
\quad
x,\xi\in\mathbf{R}^n,
\end{equation*}
and the remainder
$M_\alpha$
with the symbol
$m(x,\xi)\in S^\alpha$,
such that
\begin{eqnarray*}
m(x,\xi)
&=
&
m_1(x,\xi)+m_2(x,\xi),
\\ 
&
&
m_1(x,\xi)
=
0
,
{\rm \quad for\quad }
x\in B_1,
\ 
\xi \in \mathbf{R}^n,
\\
&
&
m_2(x,\xi)\in S^{-\infty}
{\rm \quad is \quad smoothing}.
\end{eqnarray*}
Later we will invert
$E$
modulo
a smoothing
$\psido$
via the parametrix construction.
We now rewrite
(\ref{localised_2})
as
\begin{eqnarray}
E(\eta_\rho u) 
=
-
{P}
\left(
\eta_{2\rho}
V(x) 
{Q}(\eta_\rho u)
\right)
&
+
&
\Phi_{\alpha -1}(\phi u)
\nonumber
\\
&
+
&
\Phi_{\beta -1} 
\left(
V(x) {Q}(\phi u)
\right)
\nonumber
\\
&
+
&
{P}
\left(
\eta_{2\rho}
V(x)
\Phi_{\gamma-1}
(\phi u)
\right)
\nonumber
\\
&
+
&
\Phi_\beta v
\nonumber
\\
&
+
&
f
\nonumber
\\
&
-
&
M_\alpha(\eta_\rho u ),
\label{localised_3}
\end{eqnarray}
where
all
functions 
and 
$\psido$
enjoy the same properties
as in 
(\ref{localised_2}).

{\bf \arabic{proof_step}.} 
\stepcounter{proof_step}
We claim that
\begin{equation}
\label{initial_changed}
\eta_\rho u \in W^{\sigma, r}.
\end{equation}
for some 
$\sigma$,
$r$
such that
\begin{equation}
\label{new_parameters}
\gamma< \sigma < \alpha  -\beta 
\quad
{\rm and}
\quad
\sigma - \frac{n}{r}
=
s - \frac{n}{p}.
\end{equation}
Of course,
(\ref{technical_assumption}),
(\ref{assumption2}), and
(\ref{new_parameters})
imply that
$1<r,p<\infty$.
Statement
(\ref{initial_changed})
will follow rather directly
from the Calderon-Zygmund estimates applied 
to 
(\ref{localised_3}).

Indeed,  define
$\nu=\alpha-\beta$.
Notice that due to
(\ref{assumption1})
$\gamma  < \nu$.
Rewrite
(\ref{localised_3})
as
\begin{eqnarray}
\label{abstract_form}
E(\eta_\rho u)
= 
P f_1 + \Phi_{\alpha -1} f_2 
&
+ 
&
\Phi_{\beta-1} f_3
\nonumber
\\
&
+
&
P f_4
\nonumber
\\
&
+
&
\Phi_\beta f_5
\nonumber
\\
&
+
&
f_6
\\
&
+
&
M_\alpha f_7
\end{eqnarray}
Now recall the mapping properties of
$\psido$ in $W^{s,p}$ spaces.
To be able to use them
we need 
$f_j\in L^{p_j}$
with
$1<p_j<\infty$.
This is easy to verify by applying Sobolev and Holder inequalities:
\begin{eqnarray*}
f_1
&\in& 
L^{p_1}
\quad
{\rm for}
\quad
\frac{1}{p_1}=\frac{1}{p} + \frac{\nu}{n}
-\frac{s}{n},
\\
f_2
&\in& 
L^{p_2}
\quad
{\rm for}
\quad
\frac{1}{p_2}=\frac{1}{p}-\frac{s}{n},
\\
f_3
&\in& 
W^{p_3}
\quad
{\rm for}
\quad
\frac{1}{p_3}=\frac{1}{p_1},
\\
f_4
&\in& 
L^{p_4}
\quad
{\rm for}
\quad
\frac{1}{p_4}=\frac{1}{p_1} - \frac{1}{n},
\\
f_5
&\in& 
L^{p_5}
\quad
{\rm for}
\quad
\frac{1}{p_5}=\frac{1}{q}
=
\frac{\nu}{n}
-\frac{s}{n}.
\end{eqnarray*}
Our crucial assumption
(\ref{assumption2})
ensures that all
$p_j$
are finite and 
$p_j\geq p_1>1$.

Next take a parametrix for
$E$
and apply it to  both sides of
(\ref{abstract_form}). Doing so,
$f_6$
and
$f_7$
can be ignored
because
$f_6\in C^\infty_0(\Omega)$, 
and
in the definition of
$M_\alpha$
we have
\begin{equation*}
\supp(\eta_\rho u)\cap
\supp m_1 =\emptyset.
\end{equation*}
The parametrix 
of
$E$
is a $\psido$
with the symbol in
$S^{-\alpha}$.
Utilising the mapping properties 
we discover that
\begin{equation*}
\eta_\rho u
=\sum_{j=1}^5
\tilde{f}_j
\
{\rm mod}
\
C^\infty_0(\Omega),
\end{equation*}
where
\begin{eqnarray*}
\tilde{f}_1
&\in& 
W^{\nu, p_1}
\quad
{\rm for}
\quad
\frac{1}{p_1}=\frac{1}{p} + \frac{\nu}{n} - \frac{s}{n},
\\
\tilde{f}_2
&\in& 
W^{s  +1, p_2}
\quad
{\rm for}
\quad
p_2>p_1,
\\
\tilde{f}_3
&\in& 
W^{\nu +1, p_3}
\quad
{\rm for}
\quad
{p_3}={p_1},
\\
\tilde{f}_4
&\in& 
W^{\nu , p_4}
\quad
{\rm for}
\quad
{p_4}>p_1,
\\
\tilde{f}_5
&\in& 
W^{\nu , p_5}
\quad
{\rm for}
\quad
p_5>p_1.
\end{eqnarray*}
Finally choose
$\sigma$
such that
$s<\sigma<\nu$ 
and
$\sigma<s+1$.
Applying Sobolev inequalities
we discover 
(\ref{initial_changed})
provided
\begin{equation*}
\frac{1}{r}
=
\frac{1}{p}
-
\frac{s}{n}
+
\frac{\sigma}{n} .
\end{equation*}
Thus
(\ref{initial_changed}),
(\ref{new_parameters})
are proved.
Clearly
(\ref{new_parameters})
implies
that
(\ref{assumption2})
holds for
$\sigma$ 
and 
$r$,
instead of
$s$ 
and
$p$.

%
%

{\bf \arabic{proof_step}.} 
\stepcounter{proof_step}
Our goal is to show that for some
$\sigma, r$
from
(\ref{new_parameters})
and for some
$\rho>0$
we can find
$\varepsilon>0$
such that
(\ref{initial_changed})
is improved to
\begin{equation}
\label{W^{s,p}_improved}
\eta_\rho u\in W^{\sigma+\varepsilon, r}.
\end{equation}
This is the main  statement  in  the proof.

If
(\ref{W^{s,p}_improved})
holds 
then we immediately 
conclude the proof of the theorem by
establishing
(\ref{smoothness_improved}).
Indeed, just repeat the arguments 
proving
(\ref{initial_changed})
using the improved regularity
(\ref{W^{s,p}_improved})
instead of
(\ref{new_parameters}). 
The result will be the improved integrability
of
$VQu$.
Now
(\ref{technical_assumption})
allows to apply
the
Sobolev inequality to derive that
\begin{equation*}
\eta_\rho u\in W^{s, p+\varepsilon'},
\quad
\varepsilon'>0.
\end{equation*}
The rest of the proof is devoted entirely to
the verification of
(\ref{W^{s,p}_improved}).

{\bf \arabic{proof_step}.} 
\stepcounter{proof_step}
To establish
(\ref{W^{s,p}_improved})
we will show that
by choosing
$\rho$
small
it is possible to find
$\varepsilon>0$
and a constant
$C>0$,
$C=C(\rho, u)$
such that
\begin{equation}
\label{main_line_1}
\left\|
P_k (\eta_\rho u) 
\right\|_r
\lesssim
\frac{C}{2^{(\sigma+ \varepsilon)k}}
\quad
{\rm for \quad all \quad  }
k\geq 1.
\end{equation}
Clearly it is enough to prove
(\ref{main_line_1})
only for large 
$k$. We will always assume that, say,
$k\geq 100$.
To economise on notations we set in the rest of the proof
\begin{equation*}
u=\eta_\rho u 
\quad
{\rm and}
\quad
\tilde{u}=\phi u,
\end{equation*}
cf.
(\ref{localised_3}).
Thus
$u,\tilde{u}\in W^{\sigma, r}$
and
$\supp u, \tilde{u}  \subset\subset\Omega$.
Moreover,
we can assume that
\begin{equation*}
\left\|
\tilde{u}
\right\|_{W^{\sigma, r}}\leq 1,
\end{equation*}
and hence
\begin{equation*}
\left\|
u
\right\|_{W^{\sigma, r}}\leq C_\rho,
\end{equation*}

As in the proof of 
(\ref{initial_changed})
we start  by applying the parametrix of
$E$
to
(\ref{localised_3}).
Then apply
the Littlewood-Paley
projection
$P_k$
and
take the
$L^r$-norm.
The resulting estimate will be
\begin{equation}
\label{main_line_2}
\|P_k u \|_r
\lesssim
\left\|
P_k
\Phi_{-\alpha+\beta}
\left(
\eta_{2\rho}
V 
{Q} u
\right)
\right\|_r
+
R,
\end{equation}
where we denoted by
$R$
the remainder
\begin{eqnarray*}
R
=
\left\|
P_k
\Phi_{-1}
\tilde{u}
\right\|_r
&
+
&
\left\|
P_k
\Phi_{-\alpha + \beta -1} 
\left(
V {Q}
\tilde{u}
\right)
\right\|_r
\nonumber
\\
&
+
&
\left\|
P_k
\Phi_{-\alpha+\beta}
\left(
V
\Phi_{\gamma-1}
\tilde u
\right)
\right\|_r
\nonumber
\\
&
+
&
\left \|
P_k
\Phi_{-\alpha+\beta} v
\right\|_r
\nonumber
\\
&
+
&
\frac{\lambda_N}{2^{Nk}},
\quad
{\rm with \ any\ }
N>0.
\end{eqnarray*}
Here 
$v\in L^q$
corresponds to
(\ref{localised_1}),
and the very last term 
corresponds to  the
$C^\infty_0(\Omega)$-remainders.
All
operators and functions
are compactly supported in
$\Omega$.

{\bf \arabic{proof_step}.} 
\stepcounter{proof_step}
Now we make the first step
towards proving
(\ref{main_line_1}).
Namely, we show that
the  bound 
\begin{equation}
\label{remainder_is_done}
R
\lesssim
\frac{C}{2^{ (\sigma+\varepsilon) k}}
\end{equation}
holds
for all
$k$
with
$C$
and
$\varepsilon$ 
as in 
(\ref{main_line_1}).
Indeed, let us treat the terms of
$R$
one-by-one.

Apply 
(\ref{commutator_estimate_rough})
to the first term and discover that
\begin{eqnarray*}
\left\|
P_k
\Phi_{-1} \tilde{u}
\right\|_r
&
\lesssim
&
2^{-k}
\|\Phi_{-1}P_{k-10\leq\cdot\leq k+10} \tilde{u}\|_r
+
2^{-Nk}
\\
&
\lesssim
&
2^{(-\sigma- 1)k}
+
2^{-Nk}
\\
&
\lesssim
&
\frac{1}{2^{(\sigma+1)k}}.
\end{eqnarray*}

To estimate the second term
observe that
Sobolev and Holder inequalities in combination with
(\ref{assumption2})
imply
\begin{equation*}
\left\|
V {Q}\tilde{u}
\right\|_t
\lesssim
1
\quad
{\rm for}
\quad
\frac{1}{t}=
\frac{1}{r}+
\frac{\alpha -\beta -\sigma}{n},
\quad
\frac{1}{t}
>
\frac{1}{r}.
\end{equation*}
Hence using Bernstein inequality and
(\ref{commutator_estimate_rough})
we derive  that
\begin{eqnarray*}
\left\|
P_k
\Phi_{-\alpha+\beta -1} 
\left(
V {Q}\tilde{u}
\right)
\right\|_r
&
\lesssim
&
2^
{
nk
\left(
\frac{1}{t}-\frac{1}{r}
\right)
}
\left\|
P_k
\Phi_{-\alpha+\beta -1} 
\left(
V {Q}\tilde{u}
\right)
\right\|_t
\\
&
\lesssim
&
2^{(\alpha -\beta -\sigma)k}
\left(
2^{(-\alpha+\beta-1)k}
+
2^{-Nk}
\right)
\\
&
\lesssim
&
\frac{1}{2^{(\sigma +1)k}}.
\end{eqnarray*}

To treat the next term notice that
by Holder inequality for all
$\rho>0$
\begin{equation*}
\left\|
\eta_{2\rho}
V
\Phi_{\gamma-1}
(\tilde{u})
\right\|_t
\lesssim
1
\quad
{\rm with}
\quad
\frac{1}{t}=
\frac{1}{r}
+
\frac{\alpha -\beta -\sigma -1}{n}.
\end{equation*}
There are two possibilities. If
$t\leq r$
then similarly to the previous argument
utilise
Bernstein inequality and
(\ref{commutator_estimate_rough}),
to discover that
\begin{eqnarray*}
\left\|
P_k
\Phi_{-\alpha+\beta} 
\left(
V \Phi_{\gamma -1}\tilde{u}
\right)
\right\|_r
&
\lesssim
&
2^{(\alpha -\beta -\sigma -1)k}
\left\|
P_k
\Phi_{-\alpha+\beta} 
\left(
V \Phi_{\gamma-1}\tilde{u}
\right)
\right\|_t
\\
&
\lesssim
&
2^{(\alpha -\beta -\sigma-1)k}
\left(
2^{(-\alpha +\beta)k}
+
2^{-Nk}
\right)
\\
&
\lesssim
&
\frac{1}{2^{(\sigma +1)k}}
.
\end{eqnarray*}
Otherwise we have
$t>r$. 
Then by  Holder inequality
\begin{equation*}
\|\eta_{2\rho}V \Phi_{\gamma - 1} \tilde{u}\|_r\lesssim 1
\end{equation*}
for any small 
$\rho>0$
because all functions have the compact support. 
We  deduce
from
(\ref{commutator_estimate_rough})
that
\begin{eqnarray*}
\left\|
P_k
\Phi_{-\alpha +\beta}
\left(
V
\Phi_{\gamma-1}
\tilde{u}
\right)
\right\|_r
&
\lesssim
&
2^{(-\alpha +\beta)k}
+
2^{-Nk}
\\
&
\lesssim
&
\frac{1}{2^{(\sigma+ \varepsilon) k}},
\end{eqnarray*}
where
$\varepsilon = \alpha -\beta -\sigma >0$
according to
(\ref{new_parameters}).

Finally we estimate
$\|P_k\Phi_{-\alpha +\beta}v\|_r$ 
with
$v\in L^q$.
There are two possibilities.
Suppose first that
$q\leq r$.
Then apply Bernstein inequality to discover
\begin{eqnarray*}
\|P_k \Phi_{-\alpha +\beta} v\|_r
&
\lesssim
&
2^{nk
\left(
\frac{\alpha -\beta -\gamma}{n}
-
\frac{1}{r}
\right)
}
\|{P}_k \Phi_{-\alpha +\beta} v\|_q
\\
&
\lesssim
&
2^{
\left(
\alpha-\beta-\gamma -\frac{n}{r}
\right)
k
}
2^{(-\alpha +\beta)k}
\\
&
\lesssim
&
\frac{1}{2^{(\sigma+\varepsilon) k}},
\end{eqnarray*}
where 
$\varepsilon = -\sigma +\gamma +n/r>0$
due to
(\ref{assumption2}).
Alternatively, 
suppose
$q>r$.
Recall the origin of
$v$
in
(\ref{localised_1}).
We apply
(\ref{commutator_estimate_rough})
and Holder inequality
to deduce that
\begin{eqnarray*}
\|P_k \Phi_{-\alpha +\beta} v\|_r
&
\lesssim
&
2^{(-\alpha+\beta)k}
\|v\|_r
\\
&
\lesssim
&
\frac{1}{2^{(\sigma+\varepsilon) k}}
\end{eqnarray*}
with
$\varepsilon = \alpha -\beta-\sigma>0$
due to
(\ref{new_parameters}).

Combining all estimates we conclude that
(\ref{remainder_is_done}) 
holds
with
\begin{equation}
\label{def_of_epsilon}
\varepsilon 
=
\min
\{ 
1, \alpha -\beta -\sigma, \gamma -\sigma +n/r
\},
\end{equation}
$\varepsilon>0$
because of
(\ref{assumption2})
and
(\ref{new_parameters}).

{\bf \arabic{proof_step}.} 
\stepcounter{proof_step}
It is left to estimate the first term in
(\ref{main_line_2}).
Applying 
(\ref{commutator_estimate})
from
Lemma~\ref{commutator_lemma}
we derive 
that
\begin{eqnarray}
\label{main_line_3}
\left\|
P_k\Phi_{-\alpha +\beta}
(
\eta_{2\rho}
V
Qu
)
\right\|_r
&
\lesssim
&
\left\|
\Phi_{-\alpha +\beta}
P_k
(
\eta_{2\rho}
V
Qu
)
\right\|_r
\nonumber
\\
&
&
+
\|
[P_k, \Phi_{-\alpha +\beta}]
u
\|_r
\nonumber
\\
&
\lesssim
&
2^{(-\alpha +\beta)k}
\|
\widetilde{P}_k
(
\eta_{2\rho}
V
Qu
)
\|_r
\nonumber
\\
&
&
+
2^{(-\sigma-1)k}
\|u\|_{W^{\sigma, r}}
\nonumber
\\
&
\lesssim
&
\sum_{j=k-4}^{k+4}
2^{(-\alpha + \beta)j}
\|
P_j( \eta_{2\rho} V Qu)
\|_r
\nonumber
\\
&
&
+
\frac{C_\rho}{2^{(\sigma +1)k}}
.
\end{eqnarray}
Thus our task is to estimate
$
2^{(-\alpha+\beta)k}\|P_k(\eta_{2\rho}V Qu) \|_r
$.
For convenience, in what follows we denote
\begin{equation*}
V= \eta_{2\rho} V,
\quad
\delta 
=
\|\eta_{2\rho}V\|_q.
\end{equation*}
By choosing 
$\rho$
small we can make $\delta$
as small as we wish.

To estimate
$P_k(VQu)$
we follow the standard product estimates technique,
see e.g.
\cite{Tao}. 
Taking into account the 
localisation of the 
Littlewood-Paley projections 
in the frequency space, we derive
as in
\cite{Labutin_1} 
that
for any
$k\in \mathbf{Z}$
\begin{eqnarray*}
P_k(VQu)
&
=
&
\sum_{i,j\in \mathbf{Z}}
P_k(
P_i V
P_j Qu
)
\nonumber
\\
&
=
&
\left\{
\sum_{i,j\in LL}
+
\sum_{i,j\in LH}
+
\sum_{i,j\in HL}
+
\sum_{i,j\in HH}
\right\}
P_k(
P_i V
P_j Qu
)
\nonumber
\\
&=&
I+II+III+IV,
\end{eqnarray*}
where
$LL$,
$LH$,
$HL$,
and 
$HH$
are the
low-low,
low-high,
high-low,
and
high-high
frequencies interaction zones
on the integer lattice:
\begin{eqnarray*}
LL
&=&
\left\{
i,j\in \mathbf{Z}
\colon
\
k-5\leq i,j\leq k+7,
\
\min\{i,j\}\leq k+5
\right\},
\\
LH
&=&
\left\{
i,j\in \mathbf{Z}
\colon
\
i< k-5,
\
k-3\leq j\leq k+3
\right\},
\\
HL
&=&
\left\{
i,j\in \mathbf{Z}
\colon
\
k-3\leq i\leq k+3,
\
j<k-5
\right\},
\\
HH
&=&
\left\{
i,j\in \mathbf{Z}
\colon
\
i,j>k+5,
\
|i-j|\leq 3
\right\}.
\end{eqnarray*}
We are going to estimate the  four terms 
separately.
Doing so we will constantly rely  the following 
consequence of
(\ref{commutator_estimate_rough}).
For any
$N\lesssim 1$
the inequality
\begin{equation}
\label{P_jQu}
\|P_j (Q u)\|_r
\lesssim
2^{\gamma j}
\sum_{i=j-10}^{j+10}
\left\|
{P}_i u
\right\|_r
+
{2^{-Nj}}
{C_\rho}
\end{equation}
holds for all
$j$.
As it was mentioned in
(\ref{main_line_1})
we always  assume that
$k\geq 100$.

{\bf \arabic{proof_step}.} 
\stepcounter{proof_step}
By the properties
of
$P_k$,
Bernstein inequality, and
(\ref{P_jQu})
\begin{eqnarray*}
\|I\|_r
&
\lesssim
&
\sum_{i,j\in LL}
\|
P_{i}V
\,
P_{j}Qu
\|_r
\\
&
\lesssim
&
\sum_{i,j\in LL}
\| P_{i}V \|_\infty
\|P_{j} Q u \|_r
\\
&
\lesssim
&
2^{nk/q}
\delta
\sum_{j=k-5}^{k+7}
\|P_{j} Q u \|_r
\\
&
\lesssim
&
\delta
2^{k(\alpha -\beta -\gamma)}
\sum_{j=k-20}^{k+20}
2^{\gamma j}
\|P_{j}u \|_r
\\
&
&
+\frac{C_\rho}{2^{Nk}}.
\end{eqnarray*}
Term
$II$
is estimated exactly the same way.
Thus by choosing
appropriate
$N\lesssim 1$
we derive
\begin{eqnarray}
\label{I+II}
2^{(-\alpha +\beta +\sigma)k}
(
\|I\|_r
+
\|II\|_r
)
&
\lesssim
&
\delta
\sum_{j=k-20}^{k+20}
2^{\sigma j}
\|P_j u\|_r
\nonumber
\\
&
&
+
\frac{C_\rho}{2^{100 k}}.
\end{eqnarray}

{\bf \arabic{proof_step}.} 
\stepcounter{proof_step}
To estimate
$III$
we distinguish two cases.
First, assume that
\begin{equation}
\label{r>=q}
r\geq q. 
\end{equation}
Apply the 
Holder inequality to derive 
\begin{eqnarray*}
\|III\|_r
&
\lesssim
&
\| 
(
P_{k-3\leq \cdot \leq k+3} V 
)
\,
(
P_{\cdot\leq 0} Q u
)
\|_r
\\
&
&
+
\sum_{j=1}^{k-5}
\| 
(
P_{k-3\leq \cdot\leq k+3}
)
\,
(
V 
P_{j} Q u
)
\|_r
\\
&
\lesssim
&
\| P_{k-3 \leq \cdot\leq k+3}V\|_r
\| P_{\cdot\leq 0} Q u\|_\infty
\\
&
&
+
\sum_{j=1}^{k-5}
\| P_{k-3\leq \cdot \leq k+3}V\|_r
\| P_{j} Q u\|_\infty
\\
&
=
&
X+Y.
\end{eqnarray*} 
From the Bernstein inequalities 
and
(\ref{P_jQu})
deduce that
\begin{eqnarray*}
X
&
\lesssim
&
2^
{
nk
\left(
\frac{1}{q}
-
\frac{1}{r}
\right)
}
\|V\|_{q}
\| P_{\cdot\leq 5} u\|_r
\\
&
\lesssim
&
2^
{
\left(
\alpha -\beta -\gamma
-\frac{n}{r}
\right)
k
}
\delta
C_\rho,
\end{eqnarray*}
and similarly
\begin{eqnarray*}
Y
&
\lesssim
&
\sum_{j=1}^{k-5}
2^
{
\left(
\alpha -\beta -\gamma
-\frac{n}{r}
\right)
k
}
\delta
\,
2^{
\frac{nj}{r}
}
\| P_{j} Q u\|_r
\\
&
\lesssim
&
\delta
\sum_{j=-10}^{k+10}
2^
{
\left(
\alpha -\beta -\gamma
-\frac{n}{r}
\right)
k
}
2^
{
\left(
\frac{n}{r}
+\gamma
\right)j
}
\|P_j u\|_r
\\
&
&
+
\sum_{j=-10}^{k+10}
2^
{
\left(
\alpha -\beta -\gamma
-\frac{n}{r}
\right)
k
}
\delta
\,
2^{
\frac{nj}{r}
}
\frac{C_\rho}{2^{Nj}}
.
\end{eqnarray*}
Consequently,
in the case of
(\ref{r>=q}),
after a proper choice of
$N$,
we can write the final estimate for
$III$
as
\begin{eqnarray}
\label{III_r>=q}
2^{(-\alpha +\beta +\sigma)k}
\|III\|_r
&
\lesssim
&
\delta 
\sum_{j=1}^{k+10}
2^{
\left(
\sigma -\gamma -\frac{n}{r}
\right)
(k-j)
}
\left(
2^{\sigma j}
\| P_{j} u\|_r
\right)
\nonumber
\\
&
&
+
2^
{
\left(
\sigma -\gamma -\frac{n}{r}
\right)
k
}
C_\rho
\end{eqnarray}
Notice that
$\sigma -\gamma-n/r<0$
due to
(\ref{assumption2}).

Next assume that
\begin{equation}
\label{r<q}
r<q.
\end{equation}
Hence
\begin{equation*}
\frac{1}{q}
+
\frac{1}{t}
=
\frac{1}{r}
\quad
{\rm and}
\quad
q,
t
>
r
.
\end{equation*}
By the Holder inequality
\begin{eqnarray*}
\|III\|_r
&
\lesssim
&
\| P_{k-3\leq \cdot\leq k+3}V 
\|_{q}
\|
P_{\cdot\leq 0} (Q u)\|_t
\\
&
&
+
\sum_{j=1}^{k-5}
\| 
P_{k-3\leq \cdot \leq k+3}V 
\|_q
\|
P_{j} Q u
\|_t
\\
&
=
&
Z+W.
\end{eqnarray*}
The Bernstein inequalities
imply  that
\begin{equation*}
Z
\lesssim
\delta 
C_\rho,
\end{equation*}
and
\begin{eqnarray*}
W
&
\lesssim
&
\delta
\sum_{j=1}^{k-5}
2^
{
\frac{nj}{q}
}
\|
P_{j} Qu
\|_r
\\
&
\lesssim
&
\delta
\sum_{j=-10}^{k+10}
2^{(\alpha -\beta -\gamma)j}
\,
2^{\gamma j}
\|
P_{j} u
\|_r
\\
&
&
+
\delta
\sum_{j=-10}^{k+10}
2^{(\alpha -\beta -\gamma)j}
\frac{C_\rho}{2^{Nj}}
.
\end{eqnarray*}
Consequently, in the case of
(\ref{r<q}),
the final estimate for
$III$
can be written as
\begin{eqnarray}
\label{III_r<q}
2^{(-\alpha +\beta +\sigma)k}
\| III \|_r
&
\lesssim
&
\delta
\sum_{j=1}^{k+10}
2^{(\sigma -\alpha +\beta)(k-j)}
\left(
2^{\sigma j}
\|
P_{j} u
\|_r
\right)
\nonumber
\\
&
&
+
2^{(-\alpha +\beta +\sigma)k}
C_\rho
.
\end{eqnarray}
Notice that
$-\alpha +\beta+\sigma<0$
due to
(\ref{new_parameters}).

{\bf \arabic{proof_step}.} 
\stepcounter{proof_step}
To estimate
$IV$
we 
also need to 
consider two cases.
First 
assume that
\begin{equation}
\label{r>=q'}
r\geq q'
\Leftrightarrow
\frac{1}{r}+\frac{1}{q}\leq 1.
\end{equation}
Then
define
$t\geq 1$
by writing
\begin{equation*}
\frac{1}{t}= \frac{1}{r}+ \frac{1}{q}.
\end{equation*}
Bernstein
and Holder
inequalities 
imply
that
\begin{eqnarray*}
\|P_k(P_i V\, P_j u)\|_r
&
\lesssim
&
2^
{
\frac{nk}{q}
}
\|P_k(P_i V P_j Q u)\|_t
\\
&
\lesssim
&
2^
{
(\alpha -\beta -\gamma)k
}
\delta
\|P_j Q u\|_r.
\end{eqnarray*}
Utilising
(\ref{P_jQu})
and
summing  over 
$i$
and 
$j$
in the
$HH$
region, 
we conclude 
that
in the case of
(\ref{r>=q'})
\begin{eqnarray}
\label{IV_r>=q'}
2^{(-\alpha +\beta +\sigma)k}
\|IV\|_r
&
\lesssim
&
\delta
\sum_{j=k-20}^\infty
2^{(\sigma -\gamma)(k-j)}
\left(
2^{\sigma j}\|P_j u\|_r
\right)
\nonumber
\\
&
&
+
\delta
\sum_{j=k-20}^\infty
2^{(\sigma -\gamma)k}
\frac{C_\rho}{2^{Nj}}
\nonumber
\\
&
\lesssim
&
\delta
\sum_{j=k-20}^\infty
2^{(\sigma -\gamma)(k-j)}
\left(
2^{\sigma j}\|P_j u\|_r
\right)
\nonumber
\\
&
&
+
\frac{C_\rho}{2^{100k}}.
\end{eqnarray}

Next assume that
\begin{equation}
\label{r<q'}
r<q'.
\end{equation}
By  the Bernstein and Holder inequalities
\begin{eqnarray*}
\|P_k(P_i V\, P_j  Q u)\|_r
&
\lesssim
&
2^
{
nk\left( 1-\frac{1}{r}\right)
}
\|P_k(P_i V\, P_j Q u)\|_1
\\
&
\lesssim
&
2^
{
nk\left( 1-\frac{1}{r}\right)
}
\|P_i V\|_{q}
\| P_j Q u\|_{q'}
\\
&
\lesssim
&
2^
{
nk\left( 1-\frac{1}{r}\right)
}
\delta
2^
{
nj
\left(\frac{1}{r}  -1 +\frac{1}{q}\right)}
\| P_j Q u\|_{r}.
\end{eqnarray*}
After the summation 
over 
$i$
and
$j$
lying in
the
$HH$
zone  
and the application of
(\ref{P_jQu}),
we discover that
\begin{eqnarray}
\label{IV_r<q'}
2^
{
\left(
-\alpha +\beta +\sigma
\right)
k
}
\| IV \|_r
&
\lesssim
&
\delta 
\sum_{j=k-20}^\infty
2^
{
\left(
-\alpha +\beta +\sigma
+n-\frac{n}{r}
\right)
(k-j)
}
\left(
2^{\sigma j}
\|P_j u\|_r
\right)
\nonumber
\\
&
&
+
\delta
\sum_{j=k-20}^\infty
2^
{
\left(
-\alpha +\beta +\sigma
+n-\frac{n}{r}
\right)
(k-j)
}
2^{-\gamma j}
\frac{C_\rho}{2^{Nj}}
\nonumber
\\
&
\lesssim
&
\delta 
\sum_{j=k-20}^\infty
2^
{
\left(
-\alpha +\beta +\sigma
+n-\frac{n}{r}
\right)
(k-j)
}
\left(
2^{\sigma j}
\|P_j u\|_r
\right)
\nonumber
\\
&
&
+
\frac{C_\rho}{2^{100 k}}
\end{eqnarray}
provided
(\ref{r<q'})
holds.

{\bf \arabic{proof_step}.} 
\stepcounter{proof_step}
Now we can prove the desired estimate
(\ref{main_line_1}).
First define
$\theta$
to be the smallest
parameter occurring in
(\ref{III_r>=q}),
(\ref{III_r<q}),
(\ref{IV_r>=q'}),
and
(\ref{IV_r<q'}):
\begin{eqnarray*} 
0
<
\theta
<
\min
&
\Big\{
&
\gamma +\frac{n}{r}-\sigma,
\,
\alpha -\beta -\sigma,
\,
\sigma -\gamma,
\\
&
&
-\alpha +\beta +n\left( 1-\frac{1}{r}\right) +\sigma
\Big\}.
\end{eqnarray*}
According to
(\ref{assumption2})
and
(\ref{new_parameters})
we can find such
$\theta$. 
Decrease
$\theta$,
if necessary,  so that  we also have
\begin{equation*}
\theta<\varepsilon,
\end{equation*}
where
$\varepsilon$
is taken from
from
(\ref{def_of_epsilon}).

We continue 
(\ref{main_line_2}).
The first term in the right hand side there is estimated
using
(\ref{main_line_3}),
(\ref{I+II}),
(\ref{III_r>=q}),
(\ref{III_r<q}),
(\ref{IV_r>=q'}),
and
(\ref{IV_r<q'}).
The second term in the right hand side of
(\ref{main_line_2})
is estimated by
(\ref{remainder_is_done}).
The resulting inequality is
\begin{eqnarray}
\label{theta_decay}
2^{\sigma k}
\|P_ku \|_r
\leq
C_0
\delta
\sum_{j=0}^{\infty}
\left(
2^{\sigma j}
\|P_j u\|_r
\right)
{2^{-\theta |j-k|}}
+
\frac{C(\rho)}
{2^{\theta k}}
\end{eqnarray}
where the positive constant
$C_0$
does not depend on 
$u$
and
$\rho$.

Set 
\begin{equation*}
a_k
=
2^{\sigma k}
\|P_k u\|_r,
\quad
k=0,1,\ldots
.
\end{equation*}
We will use the following  iteration 
lemma
to bound the sequence
$\{a_k\}$.

\begin{lemmy}
\label{iter_lemma}
Let
$\epsilon>0$,
let
$\delta$
satisfy
\begin{equation*}
0< \delta < (1-2^{-\epsilon})/2,
\end{equation*} 
and let a bounded  sequence
$\{ a_k \}$
satisfy
\begin{equation*}
a_k
\leq
\frac{1}{2^{\epsilon k}}
+
\delta
\sum_{j\geq 0}
\frac{a_j}{2^{2\epsilon |k-j|}}
\quad
{\rm for}
\quad
k\geq S,
\end{equation*}
with some
$S\geq 0$.
Then 
\begin{equation*}
a_k
\leq
M
\,
\|\{ a_k\}\|_{l^\infty}
\,
\frac{1}
{2^{\epsilon  k}},
\quad
k=0,1,\ldots,
\end{equation*}
with a constant
$M\geq0$,
$M=M(\varepsilon, \delta, S) $.
\end{lemmy}
This elementary lemma about number sequences 
is proved in
\cite{Labutin_1}.
The proof there is a careful but straightforward
iteration   of the  assumptions.

Going back to
(\ref{theta_decay}) we
first take
$\epsilon =\theta/2$.
Next, find
$\rho>0$
such that 
we have
\begin{equation*}
C_0 \delta
<
(
1-2^{-\epsilon/100}
)
/2.
\end{equation*}
Then in
(\ref{theta_decay})
we can choose
$S=S(u, \rho)$
such that
\begin{equation*}
a_k
\leq
\frac{1}{2^{\epsilon k}}
+
C_0 {\delta}
\sum_{j=0}^\infty
\frac
{a_j}
{2^{2\epsilon |k-j|}}
\quad
{\rm for}
\quad
k\geq S.
\end{equation*}
Finally, observe that
\begin{equation*}
a_k
\lesssim
\|u\|_{W^{\sigma,r}}
\lesssim C_\rho
\quad
{\rm for}
\quad
k=0,1,\ldots.
\end{equation*}
All assumptions of  
Lemma~\ref{iter_lemma}
now hold and 
we derive
(\ref{main_line_1}).
\end{proof}

%
%

%
%

\end{document}